\documentclass[12pt,a4paper,reqno]{article}

\usepackage{authblk}

\usepackage{graphicx}
\usepackage{caption}
\usepackage{subcaption}
\usepackage{color}
\usepackage{csquotes}

\usepackage[pdfborder={0 0 0}]{hyperref}

\usepackage{amssymb}
\usepackage{amsmath}
\usepackage{amsthm}

\usepackage[shortlabels]{enumitem}

\newtheorem{Thm}{Theorem}%[section]
\newtheorem{Prop}[Thm]{Proposition}
\newtheorem{Lem}[Thm]{Lemma}

\theoremstyle{remark}

\theoremstyle{definition}

\newcommand\set[1]{\left\{\,#1\,\right\}}		% set
\DeclareMathOperator{\spann}{span}				% span
\DeclareMathOperator{\Sp}{Sp}					% symplectic group
\DeclareMathOperator{\tr}{tr}					% trace

\def\R{\mathbb{R}}

\newcommand{\cU}{{\mathcal U}}

\usepackage[backend=biber,style=numeric, url=false, doi=false, isbn=false, giveninits=true, date=year, maxbibnames=99]{biblatex} % kein Url,Doi,Isbn, Verwende Initialien, zeige nur Jahr und keine Monate
\DeclareNameAlias{author}{first-last}					% Initialien dann Nachname
\renewbibmacro{in:}{}									% kein In: vor Journal
\DeclareFieldFormat[article,inbook,incollection,inproceedings,patent,thesis,unpublished]{title}{\textit{#1}}	% keine Anführungszeichen um Titel und kursiv
\DeclareFieldFormat[article,inbook,incollection,inproceedings,patent,thesis,unpublished]{journaltitle}{#1}
\DeclareFieldFormat{pages}{#1}							% kein pp. vor Seitenzahl
\DeclareFieldFormat{titlecase}{#1}						% keine übermäßig vielen Großbuchstaben im Titel
\DeclareFieldFormat[article]{volume}{\mkbibbold{#1}} 	% Volume in bold
\DeclareFieldFormat[article]{number}{\mkbibbold{#1}} 	% number in bold
\DeclareFieldFormat[book]{pagetotal}{}
\begin{document}

\title{A note concerning a property of symplectic matrices}
\author{Bj\"orn Gebhard\footnote{Supported by DAAD grant 57314604}}
\affil{\footnotesize Justus-Liebig-Universit\"at Gie\ss en,\\ Mathematisches Institut, Arndstr. 2, 35392 Gie\ss en}
\date{}
\maketitle

\begin{abstract}
This note provides a counterexample to a proposition stated in \cite{boscaggin_periodic_2016} regarding the neighborhood of certain $4\times 4$ symplectic matrices. 
\end{abstract}



\section{Introduction}\label{sec:intro}
%\vspace{30pt}
%\noindent
We denote by $I_n$ the $n\times n$ identity matrix, by $J$ the standard $4\times 4$ symplectic matrix, i.e.  
\[
J=\begin{pmatrix}
0 & I_2 \\
-I_2 & 0
\end{pmatrix}\in\R^{4\times 4},
\]
and by $\Sp(\R^4)=\set{S\in\R^{4\times 4}:S^TJS=J}$ the corresponding symplectic group, which shall be equipped with some norm. Furthermore a matrix $S\in\Sp(\R^4)$ is called \emph{elliptic} if the spectrum $\sigma(S)$ is contained in $S^1\setminus\{\pm 1\}$. 

In section \ref{sec:family} we will present a continuous family of symplectic matrices contradicting the following statement:
\vspace{-15pt}
\begin{displayquote}
\begin{Prop}[Prop. 2.1 of \cite{boscaggin_periodic_2016}]\label{prop:prop}
Assume that $P$ is a matrix satisfying

\vspace{-20pt}
\begin{equation}\label{eq:cond_1}
P\in\Sp(\R^4),\quad P\neq I_4,\quad \dim\ker(P-I_4)\neq 2,\quad\sigma(P)=\{1\}.
\end{equation}
Then there exists a neighborhood $\cU\subset \Sp(\R^4)$ of $P$ such that a matrix $S\in\cU$ is elliptic if and only if the following conditions hold
\begin{equation}\label{eq:cond_2}
\det(S-I_4)>0\quad\text{and}\quad \tr S < 4.
\end{equation}
\end{Prop}
\end{displayquote}
This proposition has been used in the proof of Theorem 1.1 of \cite{boscaggin_periodic_2016} to obtain a spectral stability result for periodic solutions of a perturbed Kepler problem. It has \emph{not} been used for the instability result contained in the same theorem.

\section{A family of symplectic matrices}\label{sec:family}
For $\varepsilon\geq 0$ we define the $2\times 2$ matrices
\[
A_\varepsilon=\begin{pmatrix}
1 & -\varepsilon \\
\varepsilon & 1
\end{pmatrix},\quad
B_\varepsilon=\begin{pmatrix}
1 & -\varepsilon\\
0 & 0
\end{pmatrix},\quad C_\varepsilon=\frac{1}{1+\varepsilon^2}\begin{pmatrix}
1 & -\varepsilon \\
\varepsilon & 1
\end{pmatrix}
\]
as well as the $4\times 4$ matrix
\[
P_\varepsilon=\begin{pmatrix}
A_\varepsilon & B_\varepsilon \\
0 & C_\varepsilon
\end{pmatrix}.
\]
Clearly $[0,\infty)\ni\varepsilon\mapsto P_\varepsilon\in\R^{4\times 4}$ is continuous. Moreover $P_\varepsilon \in\Sp(\R^4)$, since
\[
P_\varepsilon^TJP_\varepsilon=\begin{pmatrix}
0 & A_\varepsilon^TC_\varepsilon \\
-C_\varepsilon^TA_\varepsilon & B_\varepsilon^TC_\varepsilon-C_\varepsilon^TB_\varepsilon\end{pmatrix}=\begin{pmatrix}
0 & I_2\\
-I_2 & 0
\end{pmatrix}=J.
\]
It follows that 
\[
P_0 = \begin{pmatrix}
1 & 0 & 1 & 0 \\
0 & 1 & 0 & 0 \\
0 & 0 & 1 & 0 \\
0 & 0 & 0 & 1
\end{pmatrix}
\]
satisfies condition \eqref{eq:cond_1}. Next we show that for positive $\varepsilon$ the matrix $P_\varepsilon$ satisfies condition \eqref{eq:cond_2} and is not elliptic.
We have 
\[
\tr(P_\varepsilon)=2+\frac{2}{1+\varepsilon^2}<4.
\]
The characteristic polynomial $\chi_\varepsilon$ of $P_\varepsilon$ is given by
\[
\chi_\varepsilon(\lambda)=\left((1-\lambda)^2+\varepsilon^2\right)\left(\left(\frac{1}{1+\varepsilon^2}-\lambda\right)^2+\frac{\varepsilon^2}{(1+\varepsilon^2)^2}\right),
\]
so especially 
\[
\det(P_\varepsilon-I_4)=\chi_\varepsilon(1)=\frac{\varepsilon^4}{1+\varepsilon^2}>0
\]
and we also see that the spectrum
\[
\sigma(P_\varepsilon)=\set{1\pm i\varepsilon,(1\pm i\varepsilon)^{-1}}
\]
is not contained in $S^1$. 

Thus there exists no neighborhood of $P_0$, on which condition \eqref{eq:cond_2} implies ellipticity.

\section{Invariant Lagrangian splittings}
A Lagrangian splitting of $\R^4$ is a decomposition $\R^4=U\oplus V$ into two-dimensional subspaces satisfying 
$u_1^TJu_2=0$, $v_1^TJv_2=0$ for all $u_1,u_2\in U$, $v_1,v_2\in V$.

For $\varepsilon>0$ the planes $U=\spann\{e_1,e_2\}$, $V_\varepsilon=\spann\{v_1^\varepsilon,v_2^\varepsilon\}$, where
\[
e_1=\begin{pmatrix}
1\\0\\0\\0
\end{pmatrix},\quad
e_2=\begin{pmatrix}
0\\1\\0\\0
\end{pmatrix},\quad
v_1^\varepsilon=\begin{pmatrix}
1+\varepsilon^2\\-\varepsilon(1+\varepsilon^2)\\-2\varepsilon^2\\\varepsilon^3
\end{pmatrix},\quad
v_2^\varepsilon=\begin{pmatrix}
0\\1+\varepsilon^2\\-\varepsilon^3\\-2\varepsilon^2
\end{pmatrix},
\]
form a Lagrangian splitting of $\R^4$. Moreover, since 
\begin{align*}
P_\varepsilon e_1&=e_1+\varepsilon e_2,& P_\varepsilon e_2&=-\varepsilon e_1 + e_2,\\
 P_\varepsilon v_1^\varepsilon&=\frac{1}{1+\varepsilon^2}v_1^\varepsilon+\frac{\varepsilon}{1+\varepsilon^2}v_2^\varepsilon,& P_\varepsilon v_2^\varepsilon&=-\frac{\varepsilon}{1+\varepsilon^2}v_1^\varepsilon+\frac{1}{1+\varepsilon^2}v_2^\varepsilon,
\end{align*}
the splitting $\R^4=U\oplus V_\varepsilon$ is invariant under $P_\varepsilon$. 

On the other hand in the limiting case $\varepsilon=0$ the map $P_0$ does not admit an invariant Lagrangian splitting: Indeed let $U_0\oplus V_0$ be a splitting of $\R^4$ into $P_0$-invariant planes. We can assume that $U_0$ (otherwise $V_0$) contains a vector of the form $u_1=ae_1+be_2+e_3+ce_4$ with $a,b,c\in\R$. By the invariance also $u_2=P_0u_1=(a+1)e_1+be_2+e_3+ce_4\in U_0$. So $u_1^TJu_2=-1$ implies that the splitting is not Lagrangian.

This elaboration shows that the family $(P_\varepsilon)_{\varepsilon\in[0,\infty)}$ contradicts also a lemma on which the proof of Proposition \ref{prop:prop} is based:
\vspace{-15pt}
\begin{displayquote}
\begin{Lem}[Lem. 2.5 of \cite{boscaggin_periodic_2016}]\label{lem:lem}
Let $\{S_n\}$ be a sequence of matrices in $\Sp(\R^4)$ converging to $S$. In addition assume that for each $n\geq 0$ there exists a splitting of $\R^4$ by Lagrangian planes that are invariant under $S_n$. Then there exists another splitting by Lagrangian planes that are invariant under $S$.
\end{Lem}
\end{displayquote}

\vspace{20pt}
\noindent\textbf{Acknowledgements.} I would like to thank Prof. Rafael Ortega for his great hospitality and interesting discussions on Hamiltonian dynamics during my stay in Granada.

\printbibliography
\vspace{40pt}
%\noindent Bj\"orn Gebhard\\
% Mathematisches Institut\\
% Universit\"at Gie\ss en\\
% Arndtstr.\ 2\\
% 35392 Gie\ss en\\
% Germany\\
\noindent E-Mail address: Bjoern.Gebhard@math.uni-giessen.de
\end{document}